\documentclass[10pt]{stijl}     
  
\usepackage{amssymb} 
\usepackage{amsmath} 
\usepackage{stmaryrd}
\usepackage{euscript} 
\usepackage{latexsym} 
\usepackage{mathrsfs} 
\usepackage{wasysym}
\usepackage{float} 
\usepackage{xcolor} 
\usepackage{enumitem}
\usepackage{proof}
\usepackage{bussproofs} \EnableBpAbbreviations
\usepackage{url}
\usepackage[round]{natbib} 
\newcommand{\IPC}{{\sf IPC}}

\newcommand{\PLL}{{\sf PLL}}

\newcommand{\Kf}{{\sf K4}}

\newcommand{\Sf}{{\sf S4}}
\newcommand{\GL}{{\sf GL}}
\newcommand{\iK}{{\sf iK}_\Box}
\newcommand{\iKf}{{\sf iK4}_\Box}
\newcommand{\iKD}{{\sf iKD}_\Box}

\newcommand{\iGL}{{\sf iGL}_\Box}
\newcommand{\iSf}{{\sf iS4}_\Box}

\newcommand{\GLL}{{\sf GLL}}
\newcommand{\GPLL}{{\sf GPLL}}

\newcommand{\LJ}{{\sf G3ip}}

\newcommand{\LJX}{{\sf G3iX}}
\newcommand{\LJLL}{{\sf G3iLL}}

\newcommand{\DY}{{\sf G4ip}}

\newcommand{\DYX}{{\sf G4iX}}
\newcommand{\DYLL}{{\sf G4iLL}}

\newcommand{\DYK}{{\sf G4iK}_\Box}

\newcommand{\G}{{\sf G}}

\newcommand{\rsch}{{\EuScript R}}

\newcommand{\lgc}{{\sf L}}
\newcommand{\lang}{\ensuremath {{\EuScript L}}}
\newcommand{\defn}{\equiv _{\mbox{\em \tiny df}}} 
\newcommand{\af}{\vdash}

\newcommand{\imp}{\rightarrow} 
\newcommand{\Imp}{\ \Rightarrow\ }
\newcommand{\Ifff}{\ \Leftrightarrow\ }
\newcommand{\en}{\wedge} 
\newcommand{\of}{\vee} 

\newcommand{\E}{\exists}
\newcommand{\A}{\forall} 
\newcommand{\mdl}{\raisebox{.1mm}{$\ocircle$}} %ocircle fullmoon Circle

\newcommand{\bof}{\bigvee}
\newcommand{\ben}{\bigwedge}

\newcommand{\Ap}{\forall \hspace{-.1mm}p\hspace{.2mm}} 
\newcommand{\Ep}{\exists \hspace{-.1mm}p\hspace{.2mm}}
\newcommand{\ApR}{\forall^{\raisebox{.55mm}{\scriptsize $R$}}\hspace{-2.4mm}p\hspace{.3mm}} 
\newcommand{\EpR}{\exists^{\raisebox{.55mm}{\scriptsize $R$}}\hspace{-2.4mm}p\hspace{.3mm}}

\newcommand{\ApRi}{\forall^{\raisebox{.55mm}{\scriptsize $R^i$}}\hspace{-3.2mm}p\hspace{.3mm}} 
\newcommand{\EpRi}{\exists^{\raisebox{.55mm}{\scriptsize $R^i$}}\hspace{-3.2mm}p\hspace{.3mm}}  
\newcommand{\ApRSnot}{\forall^{\raisebox{.55mm}{\scriptsize $\overline{\rsch}$}}\hspace{-2.4mm}p\hspace{.3mm}}
\newcommand{\EpRSnot}{\exists^{\raisebox{.55mm}{\scriptsize $\overline{\rsch}$}}\hspace{-2.4mm}p\hspace{.3mm}}

\newcommand{\Appos}{\forall^{+} \hspace{-.1mm}p\hspace{.2mm}}
\newcommand{\Eppos}{\exists^+ \hspace{-.1mm}p\hspace{.2mm}}
\newcommand{\Apneg}{\forall^- \hspace{-.1mm}p\hspace{.2mm}}
\newcommand{\Epneg}{\exists^- \hspace{-.1mm}p\hspace{.2mm}}
\newcommand{\Apat}{\forall^{\it at} \hspace{-.1mm}p\hspace{.2mm}}
\newcommand{\Epat}{\exists^{\it at} \hspace{-.1mm}p\hspace{.2mm}}
\newcommand{\Apf}[1]{\forall^{\raisebox{.57mm}{\scriptsize $#1$}}\hspace{-3.6mm}p\hspace{.9mm}}
\newcommand{\Epf}[1]{\exists^{\raisebox{.57mm}{\scriptsize $#1$}}\hspace{-3.6mm}p\hspace{.9mm}}
\newcommand{\Apall}{\forall\hspace{-.1mm}p_1\dots p_n}
\newcommand{\Epall}{\exists\hspace{-.1mm}p_1\dots p_n}    

\newcommand{\seq}{\Rightarrow}
\newcommand{\sml}{\ll}
\newcommand{\smll}{\prec}

\newcommand{\intp}{\rightarrowtail}

\newcommand{\rnk}{\smll}

\newcommand{\De}{\Delta}
\newcommand{\Ga}{\Gamma}

\newcommand{\Sig}{\Sigma}

\renewcommand{\phi}{\varphi}
\newcommand{\varchi}{\raisebox{.4ex}{\mbox{$\chi$}}}

\newcommand{\cald}{{\EuScript D}}
\newcommand{\form}{{\EuScript F}}

\newcommand{\atoms}{{\cal P}}
\newcommand{\formex}{\form_{\sf ex}}

\newcommand{\MR}{{\cal R}}
\newcommand{\Gins}{\G_{\sf ins}}

\newcommand{\sqsRp}{\mathfrak{I}_{\! R}^p}
\newcommand{\sqsparRp}{\mathfrak{D}_{\! R}^p}

\newcommand{\qfset}[1]{{\sf qf}(#1)}
\newcommand{\ipa}[1]{\iota{#1}}

\newcommand{\itm}{\item[$\circ$]}

\newtheorem{Theor}{Theorem}
\newenvironment{theorem}{\begin{Theor}\em }{\end{Theor}}
\newtheorem{Lemma}{Lemma}
\newenvironment{lemma}{\begin{Lemma}\em }{\end{Lemma}}
\newtheorem{Coro}{Corollary}
\newenvironment{corollary}{\begin{Coro}\em }{\end{Coro}}
\newtheorem{Remark}{Remark}
\newenvironment{remark}{\begin{Remark}\em }{\end{Remark}}
\newtheorem{Claim}{Claim} 

\newtheorem{defin}{Definition}

\newtheorem{exam}{Example}
\newenvironment{example}{\begin{exam}\em }{\end{exam}}
\newtheorem{fct}{Fact}
\newenvironment{fact}{\begin{fct}\em }{\end{fct}}
\newenvironment{proof}{{\bf Proof}}{\hfill $\slot$}
\newcommand{\slot}{\hfill \mbox{$\dashv$}}

\newcommand{\todo}[1]{\textcolor{blue}{#1}}

\numberwithin{figure}{section}

\begin{document}
\title{Proof Theory for Lax Logic}
\vskip5pt 
\author{
Rosalie Iemhoff
\footnote{Utrecht University, the Netherlands, r.iemhoff@uu.nl. Support by the Netherlands Organisation for Scientific Research under grant 639.073.807 as well as partial support by the MOSAIC project (EU H2020-MSCA-RISE-2020 Project 101007627) is gratefully acknowledged.} }
\date{}
   
\maketitle

\todo{A mistake has been discovered in this paper, implying that the main theorem stating that \PLL\ has uniform interpolation cannot be considered proven (see Corollary 1 for the details of the mistake). A proof of the opposite is not known either. I am still trying to prove that \PLL\ has uniform interpolation, as I somehow think that is true. We'll see. To be continued.}

\begin{abstract}
 \noindent  
In this paper some proof theory for propositional Lax Logic is developed. A cut free terminating sequent calculus is introduced for the logic, and based on that calculus it is shown that the logic has uniform interpolation. Furthermore, a  separate, simple proof of interpolation is provided that also uses the sequent calculus. From the literature it is known that Lax Logic has interpolation, but all known proofs use models rather than proof systems. 
\end{abstract}



\section{Introduction}
This paper develops proof theory for an intuitionistic modal logic, Propositional Lax Logic, that models several phenomena in logic and computer science. 
When I was a PhD student, Dick taught me about intuitionistic and other intermediate logics and the numerous methods and techniques that exist to study them. Dick's love for this topic is something I have inherited from him. Many other topics to love and study would follow, but I am still grateful for Dick introducing me to this particular one. 

Propositional Lax Logic (\PLL) is an intuitionistic modal logic with a single unary modal operator $\mdl$ and is axiomatized by the principles of intuitionistic propositional logic \IPC\ plus the following three modal principles: 
\[
 \phi \imp \mdl\phi \ \ (\mdl R) \ \ \ \ \mdl\mdl\phi\imp\mdl\phi \ \ (\mdl M) \ \ \ \ 
 \mdl\phi \en\mdl\psi \imp \mdl(\phi\en\psi)\ \ (\mdl S). 
\]
In type theory, the modal operator describes a certain type constructor, and there exists a Curry-Howard like correspondence between Lax Logic and the computational lambda calculus   \citep{bentonetal1993}. In algebraic logic the operator is a so-called nucleus, and in recent work in that area, embeddings of superintuitionistic logics into extensions of \PLL\ have been used to provide new characterizations of subframe superintuitionistic logics \citep{bezhanishvilietal2019}, while in \citep{melzer2020} the method of canonical formulas for Lax Logic has been developed. 

The origins of the modality can be traced back to \cite{curry57} but the logic received its name {\em Lax Logic} in the PhD-thesis of \cite{mendler93}, where the logic was introduced in the setting of hardware verification. In that context Lax Logic is used to model reasoning about statements of the form {\em $\phi$ holds under some constraint}, denoted as $\mdl\phi$. Further investigation led to \citep{fairtlough&mendler94,fairtlough&mendler97}, in which, among other things, a sequent calculus for the logic is provided, a fact that is relevant for this paper.

In this paper we introduce a terminating sequent calculus \DYLL\ for \PLL\ which is weakly analytic (in a way explained in Section~\ref{seccalculi}) and use it to prove that the logic has uniform interpolation. 
Results on uniform interpolation for intuitionistic modal logics are rare, in particular for those that are not normal. In \citep{iemhoff17} it is shown that the logics $\iK$ and $\iKD$ have uniform interpolation, which is generalized in \citep{jalali&tabatabai2018b} to the {\sf K} and {\sf KD} versions of substructural logics weaker than \IPC, but whether other well-known intuitionistic modal logics, such as $\iKf$, $\iSf$ and $\iGL$, have uniform interpolation has, as far as we know, not been established, and the same holds for most other intuitionistic modal logics. This is in sharp contrast with the situation for classical modal and intermediate logic, about which much more is known. For example, for $\Kf$, $\Sf$ and $\GL$, which can be seen as classical counterparts of the logics mentioned above, it is known that the last one does have uniform interpolation \citep{shavrukov93} while the first two do not \citep{bilkova06,ghilardi&zawadowski95}.
We hope that this paper's proof-theoretic results for Lax Logic may lead to analogous results for other interesting members of the class of intuitionistic modal logics. 

Besides its use in the proof of uniform interpolation, the terminating sequent calculus that we develop is also interesting in connection with a paper of \cite{howe01} in which a decision method for \PLL\ is provided that is based on an annotated sequent calculus for the logic. In that paper it is shown that the calculus is sound and complete for Lax Logic and that backwards proof search in the calculus is terminating. Our sequent calculus \DYLL\ is also terminating, but consists of standard sequents instead of annotated ones. 

Finally, this work is part of the larger project of {\em Universal Proof Theory}, which aim is to develop general methods to answer universal questions about proof systems. Questions such as: When  does a theory has a certain type of proof system? What is the connection between the logical properties of a theory and its proof systems? How can it be established that a theory has no proof systems of a certain type? One possible way to obtain results in this direction is by applying a method that uses the fact that certain well-known logical properties are rare among the class of logics considered \citep{iemhoff17,jalali&tabatabai2018b}. In intermediate and (intuitionistic) modal logics, uniform interpolation is such a property. If a logic does not have uniform interpolation, then the method identifies a class of proof systems that are excluded for the logic. If, on the other hand, a logic has uniform interpolation, then that particular method does not apply. Thus we have found that the positive result that Lax Logic has uniform interpolation turns out to be a negative result in terms of the project just described.

\section{Preliminaries}
 \label{secpreliminaries}
The language $\lang$ for propositional lax logic \PLL\ contains a {\em constant} $\bot$, {\em propositional variables} or {\em atoms} $p,q,r,\dots$, the {\em modal operator} $\mdl$ and the {\em connectives} $\en,\of,\imp$, where $\neg\phi$ is defined as $(\phi \imp\bot)$. $\bot$ is by definition not an atom. A {\em circle formula} is a formula of the form $\mdl\phi$. The set of formulas in $\lang$ is denoted by $\form$. 

We denote finite multisets of formulas by $\Ga,\Pi,\De,\Sig$. $\Ga \cup \Pi$ is the multiset that contains only formulas $\phi$ that belong to $\Ga$ or $\Pi$ and the number of occurrences of $\phi$ in $\Ga \cup \Pi$ is the sum of the occurrences of $\phi$ in $\Ga$ and in $\Pi$. $\Ga,\Pi$ is short for $\Ga\cup\Pi$. Multiset $\Ga,\phi^n,\Pi$ denotes the union of $\Ga$, $\Pi$ and the multiset that consists of $n$ copies of $\phi $. We let $(\Ga \seq \De)\subseteq (\Ga' \seq \De')$ denote that $\Ga\subseteq \Ga'$ and $\De\subseteq \De'$. Furthermore ($a$ for antecedent, $s$ for succedent):  
\[
 (\Ga \seq \De)^a \defn \Ga \ \ \ \ (\Ga \seq \De)^s \defn \De \ \ \ \ 
 \mdl\Ga \defn \{ \mdl\phi \mid \phi \in \Ga \}.
\]
When working with sequents with a superscript, such as $S^i$, then $S^{ia}$ is short for $(S^i)^a$, and similarly for $S^{is}$.
We only consider (single-conclusion) sequents, which are expressions $(\Ga \seq \De)$, where $\Ga$ and $\De$ are finite multisets of formulas and $\De$ contains at most one formula. They are interpreted as $I(\Ga \seq \De) \defn (\ben\Ga \imp \bof\De)$, where $\bof\varnothing$ is interpreted as $\bot$. Expression $(\Ga \seq \De)\cdot(\Pi\seq\Sig)$ denotes $(\Ga,\Pi\seq\Sig,\De)$.
When sequents are used in the setting of formulas, we often write $S$ for $I(S)$, such as in $\af \bof_i S_i$, which thus denotes $\af \bof_i I(S_i)$. 

The {\em degree} of a formula $\phi$ is inductively defined as $d(\bot)=0$, $d(p)=1$, $d(\mdl \phi)=d(\phi)+1$, and $d(\phi\circ \psi)= d(\phi)+d(\psi)+1$ for $\circ \in\{\en,\of,\imp\}$.
In the setting of $\DYLL$ we need an order on sequents based on the {\em weight function} $w(\cdot)$ on formulas from \citep{dyckhoff92}, which is inductively defined as: the weight of an atom and the constant $\bot$ is 1, $w(\mdl \phi)=w(\phi)+1$, and $w(\phi \circ \psi) = w(\phi)+w(\psi)+i$, where $i=1$ in case $\circ \in \{\of,\imp\}$ and $i=2$ otherwise. 
We use the following ordering on sequents: $S_0 \sml S_1$ if and only if 
$S_0^a\cup S_0^s \sml S_1^a\cup S_1^s$, where $\sml$ is the order on multisets determined by weight as in \citep{dershowitz&manna79} (where they in fact define $\gg$): for multisets $\Ga,\De$ we have $\De \sml \Ga$ if $\De$ is the result of replacing one or more formulas in $\Ga$ by zero or more formulas of lower weight.

\subsection{The calculi $\LJLL$ and $\DYLL$}
 \label{seccalculi}
One of the standard calculi without structural rules for $\IPC$ is $\LJ$ (Figure~\ref{figlj}), which is the propositional part of the calculus {\sf G3i} from \citep{troelstra&schwichtenberg96}. The calculus \DY\ (Figure~\ref{figdy}) for \IPC\ from \citep{dyckhoff92} is similar to \LJ\ except for the left implication rule, which is replaced by rules that make \DY\ terminating, a feature that \LJ\ does not possess (the definition of {\em terminating} is given in Section~\ref{sectermination}). The calculus $\LJLL$ is \LJ\ extended by the two modal rules $R\mdl$ and $L\mdl$ given in Figure~\ref{figll}, which are the two modal rules of the calculus $\GPLL$ in \citep{fairtlough&mendler97}. 
The calculus $\DYLL$ is \DY\ extended by the four rules given in Figure~\ref{figll}. $\LJLL$ is an analytic calculus, meaning that it has the subformula property. $\DYLL$ does not have the subformula property, but the words of \cite{dyckhoff92} (last paragraph  Section 1) about $\DY$ apply here as well: ``Nevertheless, it has it in an obvious weak sense: one can work out what formulae are able to appear in a proof of a particular end-sequent.''
 
\begin{figure}[t]
 \centering
\[\small 
 \begin{array}{ll}
 \deduce{\Ga,p \seq p}{} \ \ \ \text{{\it Ax} \ \ ($p$ an atom)}  & 
  \deduce{\Ga,\bot\seq \De}{} \ \ \ L\bot \\
 \\ 
 \infer[R\en]{\Ga \seq \phi \en \psi}{\Ga\seq \phi & \Ga \seq \psi} & 
  \infer[L\en]{\Ga, \phi\en \psi \seq \De}{\Ga, \phi, \psi \seq \De} \\
 \\
 \infer[R\!\of \ (i=0,1)]{\Ga \seq \phi_0 \of \phi_1}{\Ga \seq \phi_i} & 
  \infer[L\of]{\Ga,\phi\of \psi\seq \De}{\Ga, \phi \seq \De & \Ga,\psi \seq \De} \\
 \\
 \infer[R\!\imp]{\Ga \seq \phi \imp \psi}{\Ga,\phi \seq \psi} & 
 \infer[L\!\imp]{\Ga,\phi\imp\psi\seq \De}{\Ga,\phi\imp\psi\seq \phi & \Ga,\psi\seq \De}\\
 \end{array}
\] 
\caption{The Gentzen calculus $\LJ$}
 \label{figlj}
\end{figure}

\begin{figure}[ht]
 \centering
\[\small 
 \begin{array}{lll}
 \deduce{\Ga,p \seq p}{} \ \ \ \text{{\it Ax} \ \ ($p$ an atom)}  & 
  \deduce{\Ga,\bot\seq \De}{} \ \ \ L\bot \\
 \\ 
 \infer[R\en]{\Ga \seq \phi \en \psi}{\Ga\seq \phi & \Ga \seq \psi} & 
  \infer[L\en]{\Ga, \phi\en \psi \seq \De}{\Ga, \phi, \psi \seq \De} \\
 \\
 \infer[R\!\of \ (i=0,1)]{\Ga \seq \phi_0 \of \phi_1}{\Ga \seq \phi_i} & 
  \infer[L\of]{\Ga,\phi\of \psi\seq \De}{\Ga, \phi \seq \De & \Ga,\psi \seq \De} \\
  \\
 \infer[R\!\imp]{\Ga \seq \phi \imp \psi}{\Ga,\phi \seq \psi} & 
 \infer[Lp\!\imp\text{ ($p$ an atom)}]{\Ga, p,p \imp \phi \seq \De}{\Ga,p,\phi \seq \De}\\
 \\
 \infer[L\en\!\imp]{\Ga, \phi\en\psi \imp \gamma \seq \De}{\Ga,\phi\imp (\psi\imp\gamma)\seq\De} & 
 \infer[L\of\!\imp]{\Ga,\phi \of \psi \imp \gamma \seq \De}{
  \Ga,\phi \imp \gamma, \psi \imp \gamma \seq \De}\\ 
 \\ \
 \infer[L\!\imp\!\imp]{Ga, (\phi\imp \psi) \imp \gamma \seq \De}{
  \Ga, \psi\imp \gamma \seq \phi \imp \psi & \gamma,\Ga \seq \De}\\
\end{array}
\] 
\caption{The Gentzen calculus $\DY$}
 \label{figdy} 
\end{figure}

\begin{figure}[h]
 \centering
\[\small 
 \begin{array}{ll}
 \infer[R\mdl]{\Ga \seq \mdl\phi}{\Ga \seq \phi} & 
 \infer[L\mdl]{\Ga,\mdl\psi\seq\mdl\phi}{\Ga,\psi\seq \mdl\phi} \\
 \\
 \infer[R\mdl^\imp]{\Ga, \mdl\phi\imp \psi \seq \De}{
  \Ga \seq \phi & \Ga, \psi \seq \De} & 
  \infer[L\mdl^\imp]{\Ga,\mdl\varchi, \mdl\phi\imp \psi \seq \De}{
  \Ga,\varchi \seq \mdl\phi & \Ga,\mdl\varchi, \psi \seq \De} 
\end{array}
\] 
\caption{$\LJLL$ is $\LJ$ plus $R\mdl$ and $L\mdl$, $\DYLL$ is $\DY$ plus all four rules.}
 \label{figll} 
\end{figure}

\subsection{Properties of calculi}
 \label{sectermination} 
A calculus is {\em terminating} with respect to an order $\smll$ on sequents if the following hold: the calculus is finite (has finitely many axioms and rules); for all sequents $S$ there are at most finitely many instances of a rule in the calculus with conclusion $S$; in every instance of a rule in the calculus the premisses come before the conclusion in the order $\smll$.  

\begin{fact}
Calculus \DYLL\ is terminating in order $\sml$ and calculus \LJLL\ is not. 
\end{fact}
 
In the section on uniform interpolation we use that \DYLL\ is a reductive calculus. Here an order on sequents $\smll$ is {\em reductive} if the following hold: it is well-founded; all proper subsequents of a sequent come before that sequent; whenever all formulas in $S$ occur boxed in $S'$, then $S \smll S'$; $(\Ga,\phi \seq \De) \smll (\Ga,q\imp\phi \seq \De)$ for all multisets $\Ga,\De$, formulas $\phi$ and atoms $q$. Clearly, $\sml$  is reductive. A calculus is {\em reductive} if it is terminating with respect to an order that is reductive. 

\begin{fact}
 \label{factreductive}
Calculus \DYLL\ is reductive with respect to order $\sml$. 
\end{fact}

A typical example of a rule that in general cannot belong to a reductive calculus is the cut rule, as in most common orders on sequents the premisses of that rule do not come before its conclusion. 

In \citep{iemhoff17} (Section 7.1), the notion of a {\em balanced calculus} \G\ is defined as: \G\ is reductive, Cut and Weakening are admissible in it, and \G\ satisfies two conditions that \DYLL\ trivially satisfies. Therefore, using Fact~\ref{factreductive} and Corollary~\ref{corequivalence} below, we can conclude the following. 

\begin{fact}
 \label{factbalanced}
Calculus \DYLL\ is a balanced calculus with respect to order $\sml$. 
\end{fact}

A rule is {\em nonflat} if its conclusion contains at least one connective or modal operator. A set of rules or a calculus is nonflat if all of its rules are. Clearly, calculi $\LJLL$ and $\DYLL$ are nonflat.

\subsection{Structural rules}
In a calculus {\em weakening is admissible} if the following two rules are admissible: 
\[
 \infer{\Ga,\phi \seq \De}{\Ga \seq \De} \ \ \ \ \ \ \ \infer{\Ga \seq \phi}{\Ga \seq \ }
\]
In a calculus {\em contraction is admissible} if the following rule is admissible: 
\[
 \infer{\Ga,\phi \seq \De}{\Ga,\phi,\phi \seq \De}
\]
A calculus has {\em cut-elimination} or {\em cut is admissible} in it if the following Cut Rule is admissible:
\[
 \infer[\it Cut]{\Ga_1,\Ga_2 \seq \De}{\Ga_1\seq \phi & \Ga_2,\phi \seq \De}
\]
In a calculus the {\em structural rules are admissible} if weakening, contraction and cut are admissible in it. 

\begin{lemma}\
 \label{lemllstructural}
Weakening and contraction as well as the following rule are admissible in $\LJLL$:
\[
 \infer{\Ga\seq\De}{\Ga\seq \bot }
\]
\end{lemma}
\begin{proof}
The proof is standard, analogous to the proof for $\LJ$ in \citep{troelstra&schwichtenberg96} and therefore omitted. 
\end{proof}

\section{Cut-elimination}
 \label{seccutelimination}
In this section we prove that the Cut Rule is admissible in $\LJLL$. The proof is by induction on the degree and a subinduction on the level of cuts. The degree of an application of the Cut Rule is the degree of the cut formula. The {\em level} of a cut in a derivation is the sum of the heights of its two premisses, where 
the {\em height} of a derivation is the length of its longest branch, where branches consisting of one node are considered to have height 1. In the rules of $\LJLL$, {\em principal formulas} are defined as usual for the rules in \LJ\ \citep{troelstra&schwichtenberg96} and for the modal rules $R\mdl$ and $L\mdl$ as given in Figure~\ref{figll}, the principal formulas are $\mdl\phi$ and $\mdl\psi$, respectively.

\begin{theorem}({\it Cut-elimination})
 \label{thmcutadmll}
\\
The Cut Rule is admissible in $\LJLL$.
\end{theorem}
\begin{proof}
Let $\af$ denote $\af_{\LJLL}$. Following the corrected version \citep{troelstra98} of the cut elimination proof for 
$\LJ$ in \citep{troelstra&schwichtenberg96}, we successively eliminate cuts from the proof, always considering those cuts that have no cuts above them, the {\em topmost} cuts. For this it suffices to show that for cut free proofs $\cald_1$ and $\cald_2$, the following proof $\cald$ of sequent $S=(\Ga_1, \Ga_2 \seq \De)$ can be transformed into a cut free proof $\cald'$ of the same endsequent: 
\[\small 
 \infer[{\it Cut}]{\Ga_1, \Ga_2 \seq \De}{
  \deduce[\cald_1]{\Ga_1 \seq \phi}{} & 
  \deduce[\cald_2]{\phi,\Ga_2 \seq \De}{} 
 }
\]
This is proved by induction on the degree of the cut formula, with a subinduction to the level of the cut. We use the fact that weakening and contraction are admissible in $\LJLL$ implicitly at various places. 

There are three possibilities:

\begin{enumerate}
\item at least one of the premises is an axiom;
\item both premises are not axioms and the cut formula is not principal in at least one of the premises;
\item the cut formula is principal in both premises, which are not axioms.  
\end{enumerate}

1.\ As in \citep{troelstra&schwichtenberg96}, straightforward, by checking all possible cases: 
If $\cald_1$ is axiom $L\bot$ or $\cald_2$ is axiom $L\bot$ and $\bot$ is not the cut formula, then we let $\cald'$ be the instance $S$ of that axiom. If $\bot$ is the cut formula, then $(\Ga_1\seq\bot)$ has a cut free derivation. Thus so has $(\Ga_1,\Ga_2\seq\De)$ by Lemma~\ref{lemllstructural}, where we use the additional fact that in this proof no new cuts are introduced with respect to $\cald_1$. Therefore we can let $\cald'$ be this proof. 

Assume both premises are not instances of $L\bot$. 
If $\cald_1$ is an instance of {\it Ax}, then $\phi$ is an atom. If $\cald_2$ also is an instance of {\it Ax}, then  $S$ is an instance of {\it Ax}, so we let $\cald'$ be that instance. If $\cald_2$ is not an axiom, then $\phi$ cannot be principal in its last inference, because it is an atom. We obtain $\cald'$ by {\em cutting at a lower level} in $\cald_2$, by which mean the following. Suppose the last inference of $\cald_2$ is $R\!\imp$:
\[
 \infer[{\it Cut}]{\Ga_1,\Ga_2\seq \psi_1\imp\psi_2}{
  \deduce{\Ga_1 \seq \phi}{\cald_1} 
  & 
  \infer[R\!\imp]{\Ga_2,\phi \seq \psi_1 \imp \psi_2}{
   \deduce{\Ga_2,\phi,\psi_1 \seq \psi_2}{\cald_2'} }
 }
\]
Consider the proof  $\cald_3$: 
\[  
 \infer[{\it Cut}]{\Ga_1,\Ga_2,\psi_1 \seq \psi_2}{
  \deduce{\Ga_1 \seq \phi}{\cald_1}
  & 
   \deduce{\Ga_2,\phi,\psi_1 \seq \psi_2}{\cald_2'} 
 }  
\]
Since the cut in $\cald_3$ is of the same degree but of lower level, there is a cut free proof $\cald_3'$ of $(\Ga_1,\Ga_2,\psi_1 \seq \psi_2)$. Hence the following derivation is the desired cut free proof $\cald'$:
\[ 
 \infer{\Ga_1,\Ga_2\seq \psi_1\imp\psi_2}{
 \deduce[\cald_3']{\Ga_1,\Ga_2,\psi_1 \seq \psi_2}{}
 } 
\]
The cases where the last inference of $\cald_2$ is another rule than $R\!\imp$ are treated in a similar way. 

The case that $\cald_1$ is not an axiom and $\cald_2$ is an instance of {\it Ax} remains. Here we also cut at a lower level, the proof is analogous to the case just treated. 

2.\ First, the case that $\phi$ is not principal in the last inference of $\cald_1$. Thus the last inference in $\cald_1$ is one of the nonmodal rules $\rsch$ of $\LJLL$ or $L\mdl$. We treat the latter, where the last part of the proof looks as follows and $\phi=\mdl\varchi$:
\[\small 
  \infer[{\it Cut}]{\Pi,\mdl\psi,\Ga_2 \seq \De}{
   \infer[L\mdl]{\Pi,\mdl\psi\seq \mdl\varchi}{
    \deduce[\cald_1']{\Pi,\psi \seq \mdl\varchi}{}
  } & 
  \deduce[\cald_2]{\Ga_2,\mdl\varchi \seq \De}{} 
  }
\] 
where $\Ga_1 = \Pi,\mdl\psi$. In case $\De$ consists of a circle formula, 
we can cut at a lower level in $\cald_1$:  
\[\small 
  \infer[L\mdl]{\Pi,\mdl\psi,\Ga_2\seq \De}{
   \infer[{\it Cut}]{\Pi,\psi,\Ga_2 \seq \De}{
    \deduce[\cald_1']{\Pi,\psi \seq \mdl\varchi}{} & \deduce[\cald_2]{\Ga_2,\mdl\varchi \seq \De}{} } 
  }
\]
Thus we obtain a proof of $(\Ga_1,\Ga_2 \seq \De)$ with a cut of the same degree as the cut in $\cald$ but of lower level, and the induction hypothesis can be applied. In case $\De$ is not a circle formula, first note that $\mdl\varchi$ cannot be principal in the last inference of $\cald_2$ as that would imply that $\De$ is a circle formula, which by assumption it is not. Thus $\mdl\varchi$ is not principal in the last inference of $\cald_2$ and 
we therefore can cut at a lower level in $\cald_2$. 

In the case that $\phi$ is not principal in the last inference of $\cald_2$, we only treat the following two cases, where $\phi=\mdl\varchi$:
\[
 \infer{\Ga_1,\Ga_2\seq \mdl\psi}{
  \deduce{\Ga_1 \seq \mdl\varchi}{} 
  &
  \infer[R\mdl]{\Ga_2,\mdl\varchi \seq \mdl\psi}{\Ga_2,\mdl\varchi \seq \psi} 
  } 
  \ \ \ \ 
 \infer{\Ga_1,\Ga_2\seq \De}{
  \deduce{\Ga_1 \seq \mdl\varchi}{} 
  &
  \infer[L\mdl]{\Pi,\mdl\varchi,\mdl\psi' \seq \De}{\Pi,\mdl\varchi,\psi' \seq \De} 
  } 
\]
where on the right side $\De$ is a circle formula and $\Ga_2=\Pi,\mdl\psi'$. Clearly, in both cases one can cut at a lower level in $\cald_2$ (cf.~case 1.).

3.\ The cut formula is principal in both premises, which are not axioms. We distinguish by cases according to the form of the cut formula and only treat the circle formulas. Thus the proof has the following form:
\[
 \infer{\Ga_1,\Ga_2\seq \mdl\psi}{
  \infer{\Ga_1 \seq \mdl\varchi}{
   \deduce{\Ga_1 \seq \varchi}{\cald_1'} } 
  &
  \infer{\Ga_2,\mdl\varchi \seq \mdl\psi}{
  \deduce{\Ga_2,\varchi\seq\mdl\psi}{\cald_2'} } 
  } 
\]
This can be replaced by a proof in which the only cut has a lower degree:
\[
 \infer{\Ga_1,\Ga_2\seq \mdl\psi}{
  \deduce{\Ga_1 \seq \varchi}{\cald_1'}  
  &
  \deduce{\Ga_2,\varchi\seq\mdl\psi}{\cald_2'}  
  } 
\]
\end{proof}

\section{Equivalence of \LJLL\ and \DYLL}
 \label{secequivalence}
In the short note \citep{iemhoff20}, a theorem has been proven that implies that \LJLL\ and \DYLL\ are equivalent. To formulate the theorem, we first need some terminology. In that paper, we define a {\em right modal rule} as a rule of the form (the $S_i$ are sequents)  
\[
 \infer[\rsch]{\Ga \seq \mdl\phi}{S_1 & \dots & S_n}
\]
Because under this definition the rule $L\mdl$ is a right modal rule, which is an unfortunate clash in terminology, we here use the name {\em succedent modal rule} for such rules instead. 
With such a succedent modal rule we associate the following implication rule
\[
  \infer[\rsch^\imp]{\Ga,\mdl\phi\imp \psi \seq \De}{
  S_1 & \dots & S_n & \Ga,\psi \seq \De}
\]
If $X$ is a set of rules, then $\LJX$ denotes the calculus $\LJ$ extended by the rules in $X$ and $\DYX$ is the calculus $\DY$ extended by the rules in $X$ {\em plus} the rules $\rsch^\imp$ for those $\rsch\in X$ that are succedent modal rules. Note that our notation \LJLL\ and \DYLL\ is consistent with this definition, that is, \DYLL\ is the result of adding the two modal rules $R\mdl$ and $L\mdl$ of \LJLL\ to \DY\ as well as the rules $R\mdl^\imp$ and $L\mdl^\imp$, which are added because both $R\mdl$ and $L\mdl$ have been added.  

In \citep{iemhoff20}, the following is shown. 

\begin{theorem} ({\it Equivalence})
 \label{thmequivalence}\\
If $\MR$ is a set of nonflat rules, the structural rules are admissible in $\LJ+\MR$, and $\DY+\MR$ is terminating and closed under weakening, then $\LJ+\MR$ and $\DYX+\MR$ are equivalent (derive exactly the same sequents) and the structural rules are admissible in $\DYX$. 
\end{theorem}

\begin{corollary}
 \label{corequivalence}
\LJLL\ and \DYLL\ are equivalent and the structural rules are admissible in both of them. In particular,  \DYLL\ is a terminating cut free sequent calculus for Lax Logic in which the Cut Rule is admissible.
\end{corollary}

\textcolor{cyan}{The last theorem is not stated correctly. In the paper \citep{iemhoff20} there is the additional requierement that $\DY+\MR$ is closed under contraction. As I do not know whether \DYLL\ is closed under contraction, I do not know whether \LJLL\ and \DYLL\ are equivalent calculi or not. This does not affect Theorem 3, that only uses \LJLL, but it undermines the proof of Theorem 6.}

\section{Interpolation}
 \label{secinterpolation} 
That \PLL\ has interpolation follows from the fact that it has uniform interpolation, as shown in the next section. Because the proof of the latter is complicated and the proof of the former is not, we present in this section a simple proof of interpolation for \PLL. Also the algorithm to compute interpolants that can be derived from it is much simpler than the one obtained from the proof of uniform interpolation. The method we present here is a straightforward  adaptation of the {\em Maehara method}, a method to prove interpolation based on sequent calculi first presented in \citep{maehara1960}, in Japanese, and later in English in \citep{ono1998}.

We introduce some terminology well-known from the literature on interpolation. Given multisets $\Ga, \De$ and formula $\phi$, let $V(\phi)$ denote all atoms that occur in $\phi$, let $V(\Ga)=\bigcup \{V(\phi) \mid \phi \in \Ga\}$, let $V(\Ga,\De)$ be short for $V(\Ga\cup\De)$, and 
$V(\Ga,\phi)$ and $V(\Ga,\phi,\psi)$ for $V(\Ga,\{\phi\})$ and $V(\Ga,\{\phi,\psi\})$, respectively. To be able to refer to the particular partition of the antecedent of a  sequent of which a formula $\varchi$ is an interpolant, we say that $\varchi$ is an {\em interpolant} of $\Ga;\Ga' \seq \De$ if $\Ga\seq \varchi$ and $\Ga',\varchi \seq \De$ are derivable in \LJLL\ and $V(\varchi)\subseteq V(\Ga)\cap V(\Ga',\De)$. 

\begin{theorem} 
 \label{thminterpolationlax}
\LJLL\ has sequent interpolation: if $\Ga,\Ga' \seq \De$ is derivable in \LJLL, then there is a formula $\varchi$ such that $\Ga\seq \varchi$ and $\Ga',\varchi \seq \De$ are derivable in \LJLL\ and all atoms in $\varchi$ occur both in $\Ga$ and in $\Ga'\cup\De$. 
\end{theorem}
\begin{proof}

We prove the theorem by constructing an algorithm that given a derivation $\cald$ of $\Ga,\Ga' \seq \De$ in \LJLL\ produces for any partition $\Ga;\Ga'\seq\De$ an interpolant. We use induction to the height of $\cald$. The case that $\cald$ is an instance of an axiom is straightforward, and so are the cases where the last inference of $\cald$ is one of the propositional rules. We treat the case that it is an implication rule as an example, followed by the case that it is one of the modal rules. 

If the last inference is $R\!\imp$, then $\De$ consists of $\phi \imp \psi$. If the partition is $\Ga;\Ga'\seq \phi\imp\psi$, then consider the partition of the premise $\Ga;\Ga',\phi \seq \psi$. By the induction hypothesis there is an interpolant $\varchi$ such that both $\Ga\seq \varchi$ and $\Ga',\varchi,\phi \seq \psi$ are derivable and $V(\varchi) \subseteq V(\Ga) \cap V(\Ga',\phi,\psi)$. This implies that $\Ga',\varchi \seq \phi\imp\psi$ is derivable and $V(\varchi)\subseteq V(\Ga',\phi \imp \psi)$. Thus proving that $\varchi$ is an interpolant for $\Ga;\Ga' \seq \phi\imp\psi$. 

If the last inference is $L\!\imp$ with principal formula $\phi \imp \psi$, then there are two possible partitions of the last sequent: 
$\Ga,\phi \imp \psi;\Ga' \seq \De$ and $\Ga;\phi \imp \psi,\Ga' \seq \De$. In the first case, by the induction hypothesis there are interpolants $\varchi_1$ and $\varchi_2$ of the partitions of the premises $\Ga';\Ga,\phi \imp \psi\seq \phi$ (note the changed order in the partition) and $\Ga,\psi;\Ga' \seq \De$. Thus 
$V(\varchi_1)\subseteq V(\Ga,\phi \imp \psi) \cap V(\Ga')$, 
$V(\varchi_2)\subseteq V(\Ga,\psi) \cap V(\Ga',\De)$, and 
the following are derivable: 
\[
 \Ga' \seq \varchi_1 \ \ \ \Ga,\phi\imp\psi, \varchi_1 \seq \phi \ \ \ \ \ \ 
 \Ga,\psi \seq \varchi_2 \ \ \ \Ga',\varchi_2  \seq \De.
\]
This implies that $\Ga', \varchi_1\imp \varchi_2 \seq \De$ and $\Ga, \phi\imp\psi,\varchi_1\seq \varchi_2$ are derivable. Thus so is $\Ga, \phi\imp\psi \seq \varchi_1\imp\varchi_2$. Since $V(\varchi_1 \imp \varchi_2)\subseteq V(\Ga,\phi \imp \psi) \cap V(\Ga',\De)$, this 
proves that $\varchi_1\imp\varchi_2$ is an interpolant of $\Ga,\phi \imp \psi;\Ga' \seq \De$.

In the second case, $\Ga;\phi \imp \psi,\Ga'\seq\De$, the induction hypothesis implies that there are interpolants $\varchi_1$ and $\varchi_2$ of the partitions of the premises $\Ga;\phi \imp \psi,\Ga'\seq \phi$ and $\Ga;\psi,\Ga' \seq \De$.
Thus $V(\varchi_1)\subseteq V(\Ga) \cap V(\Ga',\phi \imp\psi)$, 
$V(\varchi_2)\subseteq V(\Ga) \cap V(\Ga',\psi,\De)$, and the following are derivable: 
\[
 \Ga \seq \varchi_1 \ \ \ \Ga',\phi \imp \psi, \varchi_1 \seq \phi \ \ \ \ \ \ 
 \Ga \seq \varchi_2 \ \ \ \Ga', \psi,\varchi_2  \seq \De.
\]
It is not hard to see that $\varchi = (\varchi_1\en\varchi_2)$ is the desired interpolant. 

If the last inference is $R\mdl$, then for the interpolant of the conclusion one can take the interpolant of the premise.

If the last inference is $L\mdl$, then $\De$ consists of a circle formula $\mdl\psi$ and the principal formula is $\mdl\phi$. For partition $\Ga;\Ga',\mdl\phi \seq\mdl\psi$ of the conclusion, let $\varchi$ be the interpolant for the corresponding partition of the premise, namely such that the following are derivable:
\[
 \Ga \seq \varchi \ \ \ \ \Ga',\phi,\varchi \seq \mdl\psi. 
\] 
Clearly, $\varchi$ is an interpolant for the conclusion as well. 
Similarly, for partition $\Ga,\mdl\phi;\Ga'\seq\mdl\psi$ of the conclusion, let $\varchi$ be the interpolant for the corresponding partition of the premise, namely such that the following are derivable:
\[
 \Ga,\phi \seq \varchi \ \ \ \ \Ga',\varchi \seq \mdl\psi. 
\] 
If $\varchi$ is a circle formula, then it is an interpolant for the conclusion. Otherwise we can take $\mdl\varchi$ instead. 
\end{proof}

\begin{corollary} 
 \label{corinterpolationlax}
Lax Logic has Craig interpolation. In particular, there is an algorithm that given a derivation of $\phi \seq \psi$ in \LJLL\ produces an interpolant for $\phi \imp \psi$. 
\end{corollary}

\section{Uniform interpolation}
In \citep{iemhoff17} a modular and constructive method was developed to prove uniform interpolation for certain intermediate and intuitionistic modal logics, which was later generalized to other substructural (modal) logics by \cite{jalali&tabatabai2018a,jalali&tabatabai2018b}. In this section we prove that Lax Logic, a logic not covered by these previous methods, has uniform interpolation. Our proof is an adaptation of the proof in \citep{iemhoff17} that the intuitionistic modal logic $\iK$  has uniform interpolation. For background information and motivation we refer the reader to that paper. 

\subsection{Uniform interpolants for formulas}
 \label{secuniforminterpolants}
A logic has {\em uniform interpolation} if for any atom $p$ and any set of atoms 
$\atoms$ not containing $p$, 
the embedding of $\form (\atoms)$ into $\form (\atoms \cup \{p\})$ has a right and a left 
adjoint: For any formula $\phi$ and any atom $p$ there exist two formulas, usually denoted by $\Ap \phi$ and $\Ep \phi$, in the language of the logic, that do not contain $p$ and such that for all $\psi$ not containing $p$: 
\[
 \af \psi \imp \phi \Ifff\ \af \psi \imp \Ap \phi \ \ \ \ \af \phi \imp \psi \Ifff\ \af \Ep\phi \imp \psi. 
\] 
Given a formula $\phi$, its {\em universal uniform interpolant with respect to $p_1\dots p_n$} is $\Apall\phi$, which is short for $\A p_1 (\A p_2 (\dots (\A p_n \phi)\dots )$, and its {\em existential uniform interpolant with respect to $p_1\dots p_n$} is $\Epall\phi$, short for $\E p_1 (\E p_2 (\dots (\E p_n \phi)\dots )$. The requirements above could be replaced by the following four requirements.   
\[ 
 \tag{$\A$} \label{eqexu}
  \af \Ap\phi \imp \phi \ \ \ \ \af \psi \imp \phi \Imp \af \psi \imp \Ap \phi. 
\]
\[ 
 \tag{$\E$} \label{eqexe}
  \af \phi \imp \Ep \phi\ \ \ \ \af \phi \imp \psi \Imp \af \Ep\phi \imp \psi. 
\]
In classical logic one only needs one quantifier, as $\Ep$ can be defined as $\neg\Ap \neg$ and vice versa. Although in the intuitionistic setting $\Ep$ can also be defined in terms of $\Ap$, namely as $\Ep\phi = \A q(\Ap(\phi \imp q) \imp q)$ for a $q$ not in $\phi$, having it as a separate quantifier is convenient in the proof-theoretic approach presented here (we follow \citep{pitts92}, which also uses both quantifiers).  

In order to prove that Lax Logic has uniform interpolation, the first task is to adapt the above notions to the setting of sequents and sequent calculus \DYLL.

\subsection{Reductive orders and rank}
The set $\formex$ is the smallest set of expressions that contains all formulas in the language $\lang$, is closed under the connectives and modal operator, and if $S$ is a sequent in $\lang$ and $p$ an atom, then $\Ap S$ and $\Ep S$ belong to $\formex$. For example, when $S$ is a sequent in $\lang$ and $\phi$ a propositional formula, then $(\phi \imp \Ep S)$ belongs to $\formex$, as does $\mdl(\phi \en \Ap S)$, but $\Ep \E q S$ does not.  

The reductive order $\sml$ on sequents defined in Section~\ref{secpreliminaries} determines the following order on formulas in $\formex$, denoted $\smll$. For any expression $\phi$ in $\formex$, let $\qfset \phi$ denote the multiset consisting of all occurrences of subformulas of the form $QpS$ in $\phi$, where $Q \in \{\E,\A\}$. The order on multisets of the form $\qfset \phi$ again is in the style of \citep{dershowitz&manna79}: $\qfset \phi \smll_{\sf qf} \qfset \psi$ iff $\qfset \phi $ is the result of replacing one or more formulas of the form $QpS$ in $\qfset \psi$ by zero or more formulas of the form $Q'qS'$ with $S' \sml S$, where $Q,Q' \in \{\E,\A\}$. This order is well--defined since by definition such $S$ are sequents in $\lang$ and therefore can be compared via $\sml$. The order on $\formex$ can now be defined as: if $\phi,\psi \in \form$, then $\phi \smll \psi$ iff  $(\, \seq \phi) \smll (\,\seq \psi)$; if $\phi \in \form$ and $\psi\not\in\form$, then $\phi \rnk  \psi$ and not $\psi \rnk \phi$; if $\phi,\psi \not\in \form$, then $\phi \rnk \psi$ if $\qfset\phi \smll\qfset\psi$. When $\phi \rnk \psi$, we say that $\phi$ is of {\em lower rank\/} than $\psi$. Clearly, since the order $\sml$ on sequents is well--founded, so is the order $\rnk$ on $\formex$.

\subsection{Interpolant assignments}
 \label{secassignments}
For any given calculus, let $\Gins$ denote the set of instances of rules in \G. 
A sequent $S$ is {\em principal} for an instance $R$ of a rule if the conclusion of $R$ is of the form $S'\cdot S$ for some sequent $S'$ and the principal formula of $R$ belongs  to $S$. Otherwise $S$ is {\em nonprincipal}. For example, suppose $R$ has conclusion $(\Ga,\phi \seq \De)$ and $\phi$ is the principal formula of $R$, then any sequent of the form $(\Ga',\phi \seq \De')$, where $(\Ga' \seq \De')\subseteq (\Ga\seq \De)$, is principal for $R$. 

An {\em interpolant assignment} $\ipa{}$ for $\G$, assigns, for every atom $p$ and sequent $S$, 
$\ipa \Ep S=\top$ and $\ipa \Ap S=\bot$ in case $S$ is empty, and in case $S$ is not empty:
  
\begin{itemize}

\itm for every $R\in \Gins$ with conclusion $S$, to each of the expressions $\EpR S$ and $\ApR S$ a formula in $\formex$ that is of lower rank than $\Ep S$ (or, equivalently, of lower rank than $\Ap S$), which are denoted by $\ipa{\EpR S}$ and $\ipa{\ApR S}$, respectively, and 
\itm for every $\rsch \in \G$ such that $S$ is nonprincipal for at least one instance of $\rsch$, to each of the expressions $\EpRSnot S$ and $\ApRSnot S$ a formula  
in $\formex$ that is of lower rank than $\Ep S$, which are denoted by $\ipa{\EpRSnot S}$ and $\ipa{\ApRSnot S}$, respectively.

\end{itemize}

Given an interpolant assignment we define the following formulas in $\formex$. Recall that $p$ and $q$ range over atoms. 
\renewcommand*{\arraystretch}{1.3}
\[
 \begin{array}{lll}
  \Appos S & \defn & \bof \{ \ipa{\ApR S} \mid R \in \Gins, \text{ $S$ is the conclusion of $R$}\} \\
  \Apneg S & \defn & \bof \{\ipa{\ApRSnot S} \mid \text{$\rsch \in \G$, $S$ is nonprincipal for some instance of $\rsch$}\} \\
  \Eppos S & \defn & \ben\{\ipa{\EpR S} \mid R \in\Gins, \text{ $S$ is the conclusion of $R$} \} \\
  \Epneg S & \defn & \ben\{\ipa{\EpRSnot S} \mid\text{$\rsch \in \G$, $S$ is nonprincipal for some instance of $\rsch$}\} \\
  \Apat  S & \defn & \bof \{q \in S^s \mid \text{$q$ an atom and $q \neq p$, or $q=\top$} \} \, \of \\
           &       & \bof \{ q \en \Ap (\phi, S^a \backslash \{q \imp \phi\} \seq S^s) \mid 
                         (q \imp \phi) \in S^a, q \neq p \} \\
  \Epat  S & \defn & \ben \{q \in S^a \mid \text{$q$ an atom and $q \neq p$, or $q=\bot$} \} \, \en \\
           &       & \ben \{q \imp \Ep(\phi, S^a \backslash \{q \imp \phi\} \seq S^s) \mid 
                           (q \imp \phi) \in S^a, q \neq p \}. 
 \end{array}
\]
Observe that there could be more than one instance of a single rule $\rsch$ 
that has $S$ as a conclusion, in which case every instance corresponds to a separate disjunct or conjunct of the interpolant assignment. The definition above is well-defined for reductive calculi, because for such calculi all sets over which the big conjunctions and disjunctions range are finite. 

We define a rewrite relation $\intp$ on $\formex$ that is the smallest relation on $\formex$ that preserves the logical operators and satisfies:   
\[ 
 \Ap S \ \intp \ \Appos S \of \Apneg S \of \Apat S 
 \ \ \ \ 
 \Ep S \ \intp \ \Eppos S \en \Epneg S \en \Epat S. 
\]
After the definition of uniform interpolation for sequents in Section~\ref{secuipseq} below, 
three examples of the rewrite relation are provided.  

As $\DYLL$ is a reductive calculus, the following lemma follows from Lemma 3 in \citep{iemhoff17}.

\begin{lemma}
 \label{lemunqiue}
If there exists an interpolant assignment for \DYLL, then the relation $\intp$ on $\formex$ is confluent and strongly normalizing. 
\end{lemma} 

Note that according to this relation, any formula in $\lang$ is in normal form, and any expression in $\formex$ that is not a formula in $\lang$ is not.
Let $\overline{\phi}$ denote the normal form of an expression $\phi$. For any sequent $S$ in $\lang$ we now define its {\em left} and {\em right} {\em interpolant} to be $\overline{\Ap S}$ and $\overline{\Ep S}$, respectively. We omit the overline most of the time, as it will be clear from the context whether the expression in $\formex$ is meant or the actual normal form in $\lang$, i.e.\ the interpolants.

\subsection{Uniform interpolants for sequents}
 \label{secuipseq}
To define uniform interpolants in the setting of sequents, we introduce the notion of a partition, which applies to sequents and to rules. 
A {\em partition} of a sequent $S$ is an ordered pair $(S^r,S^i)$ ($i$ for {\em interpolant}, $r$ for {\em rest}) such that $S = S^r \cdot S^i$. It is a {\em $p$--partition} if $p$ does not occur in $S^r$. 
For any sequent $S$ and partition $(S^i,S^r)$ we use the abbreviation:
\[
 S^r \cdot (\Ep S^i \seq \Ap S^i \mid \varnothing) 
 \defn 
 \left\{ 
  \begin{array}{ll}
   S^r \cdot (\Ep S^i \seq \Ap S^i) & \text{if $S^s\neq\varnothing$ and $S^{rs}=\varnothing$} \\
   S^r \cdot (\Ep S^i \seq \ ) & \text{if $S^s=\varnothing$ or $S^{rs}\neq\varnothing$.}
  \end{array}
 \right.
\]
A {\em (p-)partition} of an instance $R =(S_1 \dots S_n/S_0)$ of a rule is a (p-)partition of the sequents in the rule. Given such a partition and $\star\in \{r,i\}$, let $(R^r,R^i)$ and $R^\star$  respectively denote the expressions  
\[
 \infer[(R^r,R^i)]{(S_0^r,S_0^i)}{(S_1^r,S_1^i) & \dots & (S_n^r,S_n^i)} 
 \ \ \ \ \ \ \ 
 \infer[R^\star]{S_0^\star}{S_1^\star & \dots & S_n^\star}
\]
In Lemma 2 in \citep{iemhoff17} it is shown that given an interpolant assignment an intuitionistic modal logic has uniform interpolation if the following holds:

\begin{itemize}
\item[($\A$l)] 
For all $p$: $\af S^a,\Ap S \seq S^s$; 

\item[($\E$r)] 
For all $p$: $\af S^a \seq \Ep S$; 

\item[($\A\E$)] 
If $S$ is derivable, for all $p$ and all $p$--partitions $(S^r,S^i)$: 
$\af S^r \cdot (\Ep S^i \seq \Ap S^i \mid \varnothing)$. 
\end{itemize}
 
Properties ($\A$l) and ($\E$r) are the {\em independent} (from partitions) {\em interpolant properties}, and ($\A\E$) is the {\em dependent interpolant property}. It is also shown that $\Ap S$ and $\Ep S$ are equivalent to $\Ap I(S)$ and $\Ep (\ben S^a)$, respectively. In particular,  $\Ap (\, \seq \phi)$ is equivalent to $\Ap \phi$ and $\Ep (\phi \seq \,)$ to $\Ep \phi$.

A partition $(S^r,S^i)$ of $S$ {\em satisfies} the interpolant properties if, in the case of the independent property, $S$ satisfies them (in which case we also say that $S$ satisfies them), and in case of the dependent property, it holds for that particular partition. A sequent {\em satisfies} a property if every possible partition of the sequent satisfies it. 

To show that Lax Logic has uniform interpolation it thus suffices to show that we can define, for all sequents $S$, formulas $\Ap S$ and $\Ep S$ in $\lang$ in such a way that the three interpolant properties hold. Before we do so, three examples of uniform sequent interpolants are given.

\begin{example}
Suppose the calculus only contains the rule $\rsch$ for conjunction on the left: 
\[ 
 \infer{\Ga , \phi \en \psi \seq \De}{\Ga,\phi, \psi \seq \De}
\]
We provide the uniform interpolants for $S=(p\en q,r,s\seq t)$. We have 
\[
 \Ap S \intp \iota\ApR S \of \ipa{\ApRSnot S} \of \iota\Apat S \ \ \ \ 
 \Ep S \intp \iota\EpR S \en \ipa{\EpRSnot S} \en \iota\Epat S,
\]
where $R=S_1/S$ stands for the instance of $\rsch$ with $p\en q$ as the principal formula. Thus $S_1 = (p,q,r,s\seq t)$. 
The standard interpolant assignment introduced below satisfies $\iota \ApR S = \Ap S_1$, $\iota \EpR S = \Ep S_1$, $\ipa{\ApRSnot S}=\bot$, and $\ipa{\EpRSnot S}=\top$. This implies that  
\[
 \Ap S \intp \Ap S_1\of \bot \of t \ \ \ \ \Ep S \intp \Ep S_1 \en \top \en r \en s. 
\]
Since $S_1$ cannot be the conclusion of any rule, we have $\Ap S_1=\ipa{\ApRSnot S} \of \Apat S_1 = \bot \of t$ and $\Ep S_1=\ipa{\EpRSnot S} \en \Epat S_1 = \top \en q \en r\en s$. Therefore 
\[
 \Ap S \intp \bot \of t \of \bot \of t \ \ \ \ \Ep S \intp \top \en q \en r \en s \en \top \en r \en s.  
\]
Indeed, with $\Ap S =t$ and $\Ep S = q \en r\en s$, the three interpolant properties are satisfied for $S$.  
\end{example}

\begin{example}
For the same calculus as in the previous example, the uniform interpolants of the sequent $S=(p_0\en q_0, p_1 \en q_1 \seq r)$ are computed as follows, where we abbreviate $p_i \en q_i$ with $\psi_i$. For $i=0,1$, let $R_i$ stand for the instance of $\rsch$ with $\psi_i$ as the principal formula, and define sequents $S_0 = (p_0, q_0,\psi_1 \seq r)$ and $S_1 = (\psi_0, p_1,q_1 \seq r)$, which are the premises of these instances. 
By the above definition,
\[
 \Ap S \intp \iota\Apf{R_0}S \of \iota\Apf{R_1}S \of \ipa{\ApRSnot S} \of \iota\Apat S \ \ \ \ 
  \ \ \ \ 
 \Ep S \intp \iota\Epf{R_0}S \en \iota\Epf{R_1}S \en \ipa{\EpRSnot S} \en \iota\Epat S.
\]
The disjuncts and conjuncts can be computed as in the previous example. 
\end{example}

\begin{example}
Suppose the calculus consists only of the rule $\rsch$ for disjunction on the right:
\[
 \infer[(i=0,1)]{\Ga \seq \phi_0 \of \phi_1}{\Ga \seq \phi_i}
\]
We provide the uniform interpolants for $S=(r\seq p \of q)$. Let $R_p, R_q$ stand for the instances of $\rsch$ with conclusion $S$ and premises $S_p = (r \seq p)$ and $S_q = (r\seq q)$, respectively. By the above definition,
\[
 \Ap S \intp \iota\Apf{R_p}S \of \iota\Apf{R_q}S \of \ipa{\ApRSnot S} \of \iota\Apat S
 \ \ \ \ 
 \Ep S \intp \iota\Epf{R_p}S \en \iota\Epf{R_q}S \en \ipa{\EpRSnot S} \en \iota\Epat S.
\]
The standard interpolant assignment introduced below satisfies $\iota \ApR S = \Ap S_1$ for any instance $R=S_1/S$ of $\rsch$, as well as $\ipa{\ApRSnot S}=\bot$ and $\ipa{\EpRSnot S}=\top$ for any $S$. This, together with the fact that $\iota\Apat S=\bot$, implies 
\[
 \Ap S \intp \iota\Ap S_p \of \iota\Ap S_q \of \bot \of \bot 
 \ \ \ \ 
 \Ep S \intp \iota\Ep S_p \en \iota\Ep S_1 \en \top \en r.
\]
Since $S_p$ and $S_q$ cannot be the conclusion of any rule, we have 
\begin{eqnarray*}
 \Ap S_p=\ipa{\ApRSnot S_p} \of \Apat S_p = \bot \of \bot \ \ \ \ 
  \Ap S_q=\ipa{\ApRSnot S_q} \of \Apat S_q = \bot \of q
 \\
 \Ep S_p=\ipa{\EpRSnot S_p} \en \Epat S_p = \top \en r  =\ipa{\EpRSnot S_q} \en \Epat S_q = \Ep S_q.
\end{eqnarray*}
Therefore $\Ap S \intp \bot \of \bot \of \bot \of q \of \bot \of \bot$ and  
$\Ep S \intp \top \en r \en \top \en r \en \top \en r$. Indeed, with $\Ap S =q$ and $\Ep S = r$, the three interpolant properties are satisfied for $S$.
\end{example}

\subsection{The inductive properties}
 \label{secinductive} 
For a modular development of uniform interpolation, we introduce the following six properties of rules, where $\varnothing \af \phi$ should be read as $\af\phi$. 
Given an instance $R = (S_1\dots S_n/S_0)$ of a rule $\rsch$, we define 
\begin{eqnarray*}
 \sqsRp    & \defn & \{ S_j \cdot (\Ap S_j \seq \,), (S_j^a\seq \Ep S_j) \mid 
                     1 \leq j \leq n \} \cup \\
           &       & \{ S^a \seq \Ep (S^a \seq \,) \mid 
                      S \subset S_0 \text{ or } \mdl S \subseteq S_0  
                      \text{ or } S \subseteq S_j \text{ for some }1\leq j\leq n \}  \\ 
 \sqsparRp & \defn & \bigcup_{j=1}^n \{ S_{j}^r\cdot 
                     (\Ep S_{j}^i \seq \Ap S_{j}^i \mid \varnothing) \mid 
                     \text{$(S_j^r,S_j^i)$ a $p$--partition of $S_j$} \}.
\end{eqnarray*}
For $\sqsparRp$, note that it contains the sequent $S_j^r\cdot (\Ep S_j^i \seq \Ap S_j^i \mid \varnothing)$ for {\em any} possible $p$-partition $(S_j^r,S_j^i)$ of a premiss $S_j$ of $R$. And that for $S$ with empty succedent, $S^r\cdot (\Ep S^i \seq \,)$ derives $S^r\cdot (\Ep S^i \seq \Ap S^i)$.  

The sets $\sqsRp$ and $\sqsparRp$ contain the sequents to which, in a proof of the interpolant properties that uses induction along $\smll$, the induction hypothesis could be applied. In such a proof, the assumption that the interpolant properties hold for all sequents below $S$ implies that the sequents in $\sqsRp$ and $\sqsparRp$ are derivable. Note that in the case that $R$ is an instance of an axiom, the latter set is empty, but the former is not, as it contains all sequents of the form $S^a \seq \Ep (S^a \seq \,)$ for $S$ such that $S \subset S_0$ or $\mdl S \subseteq S_0$.   
 
\begin{description}[leftmargin=4.2em,style=nextline]

\item[\rm (IPP)$_\rsch^\A$]
$\sqsRp \af S \cdot (\ApR S\seq \,)$ for every instance $R$ of $\rsch$ with conclusion $S$.  

\item[\rm (IPN)$_\rsch^\A$] 
If $S$ is nonprincipal for some instance of $\rsch$, then the assumption that all sequents below $S$ satisfy the interpolant properties implies 
$\af S \cdot (\ApRSnot S\seq \,)$.   

\item[\rm (IPP)$_\rsch^\E$]
$\sqsRp \af (S^a \seq \EpR S)$ for every instance $R$ of $\rsch$ with conclusion $S$. 

\item[\rm (IPN)$_\rsch^\E$] 
If $S$ is nonprincipal for some instance of $\rsch$, then the assumption that all sequents below $S$ satisfy the interpolant properties implies 
$\af (S^a \seq \EpRSnot S)$.  

\item[\rm (DPP)$_\rsch$]
For every sequent $S$ that has a derivation of which the last inference is an instance $R$ of $\rsch$, and for every $p$--partition $(S^r,S^i)$ such that sequent $S^i$ is principal for $R$:  
$\sqsparRp \af S^r \cdot (\Ep S^i \seq \Ap S^i \mid \varnothing)$.  

\item[\rm (DPN)$_\rsch$]
For every sequent $S$ that has a derivation of which the last inference is an instance $R$ of $\rsch$, and for every $p$--partition $(S^r,S^i)$ such that sequent $S^i$ is nonprincipal for $R$: if all sequents that are below $S$ satisfy the interpolant properties, then $\af S^r \cdot (\Ep S^i \seq \Ap S^i \mid \varnothing)$. 
\end{description}

These six properties are called the {\em inductive properties} in this paper. 
``IP'' stands for {\em independent property}, ``DP'' for {\em dependent property},  
``P'' and ``N'' for  {\em principal} and {\em nonprincipal}, respectively.

For an intuitionistic modal logic \lgc\ with calculus \G, an interpolant assignment is {\em sound} for a rule $\rsch$ in \G, if the six inductive properties hold for $\rsch$, where ``$\af$'' equals ``$\af_\lgc$''. It is {\em sound} for a calculus if it is sound for all the rules of the calculus. 

\begin{remark}
 \label{remcondis}
The following observation will be used to prove (DPP)$_\rsch$ and (DPN)$_\rsch$. Consider a sequent $S$ with partition $(S^r,S^i)$, which has a derivation of which the last inference is an instance $R=(S_1\dots S_n/S)$ of $\rsch$. To prove (DPP)$_\rsch$, thus in case $S^i$ is principal for $R$, in order to prove 
\[
 \sqsparRp \af S^r \cdot (\Ep S^i \seq \Ap S^i \mid \varnothing)
\] 
it suffices to show that 
\[
 \sqsparRp \af S^r \cdot (\EpRi S^i \seq \ApRi S^i \mid \varnothing)
\] 
for some partition $(R^r,R^i)$ of $R$ with conclusion $(S^r,S^i)$ such that $R^i$ is an instance of $\rsch$. The reason being, that for such an $R^i$, $\EpRi S^i$ is a conjunct of $\Ep S^i$ and $\ApRi S^i$ a disjunct of $\Ap S^i$. 
Likewise, to prove (DPN)$_\rsch$, thus in case $S^i$ is nonprincipal for $R$, to prove that 
\[
 \af S^r \cdot (\Ep S^i \seq \Ap S^i \mid \varnothing)
\] 
it suffices to prove that 
\[
 \af S^r \cdot (\EpRSnot S^i \seq \ApRSnot S^i \mid \varnothing).
\] 
\end{remark}

Corollary 32 and Theorem 31 in \citep{iemhoff17} state the following. 

\begin{theorem}
 \label{thmsufficientmod}
Any intuitionistic modal logic \lgc\ with a balanced calculus that consists of $\DYK$,  focused rules, and modal focused rules, has uniform interpolation.
\end{theorem}

\begin{theorem}
 \label{thmsufficientmod}
Any intuitionistic modal logic $\lgc$ with a balanced calculus that contains $\DYK$ and has an interpolant assignment that is sound with respect to all rules that neither are focused, nor  modal focused, nor belong to $\DYK$, has uniform interpolation. 
\end{theorem}

Given the fact that all rules of \LJ\ are focused, as shown in \citep{iemhoff17}, to prove that \PLL\ has uniform interpolation it suffices to show that the four extra rules of \DYLL\ are sound with respect to an interpolant assignment that is an extension of the one from \citep{iemhoff17}, which is called the {\em standard interpolant assignment}. This is shown in the next section.

\subsection{Interpolant assignment for \PLL}
 \label{secassignment}
Here the interpolant assignment for the four modal rules of $\DYLL$ is introduced. Its soundness is proved in the next section. 

For $R = S_1/S$ being an instance 
of $\rsch=R\mdl$ or $\rsch=L\mdl$, the interpolant assignment is defined as follows:
\[
 \begin{array}{lll}
  \EpR S \defn \mdl\Ep S_1 & & \ApR S\defn \mdl\Ap S_1 \\
  \EpRSnot S \defn \top    & & \ApRSnot S\defn \bot \\
 \end{array}
\]
For $R$ being an instance 
\[
  \infer[=\ ]{S}{S_1 & S_2}
  \infer{\Ga,\mdl\phi \imp \psi \seq \De}{\Ga\seq \phi & \Ga,\psi \seq \De}
\]
of $\rsch=R\mdl^\imp$, the interpolant assignment is defined as follows:
\[
 \begin{array}{llll}
  \EpR S \defn \Ep S_1 \en (\Ap S_1 \imp \Ep S_2) & & \ApR S \defn \Ap S_1 \en \Ap S_2 \\
  \EpRSnot S \defn \top    & & \ApRSnot S\defn \bot & \text{if $S^s = \varnothing$} \\
  \EpRSnot S \defn \Ep (S^a \seq \, )  & & \ApRSnot S\defn \bot & \text{if $S^s \neq \varnothing$} \\
 \end{array}
\]
For $R$ being an instance 
\[
  \infer[=\ ]{S}{S_1 & S_2}
  \infer{\Ga,\mdl\varchi,\mdl\phi \imp \psi \seq \De}{\Ga, \varchi \seq \mdl\phi & \Ga,\mdl\varchi,\psi \seq \De}
\]
of $\rsch=L\mdl^\imp$, the interpolant assignment is defined as follows:
\[
 \EpR S \defn \mdl\Ep S_1 \en (\mdl\Ap S_1 \imp \Ep S_2) \ \ \ \ 
 \ApR S \defn \mdl\Ap S_1 \en \Ap S_2.
\]
For $\EpRSnot S$ and $\ApRSnot S$, first define for $\gamma = \mdl\alpha\imp\beta$: 
\begin{eqnarray*}
 & & S^{\gamma 0} \defn (S^a\backslash \{\gamma\} \seq \mdl\alpha) \\
 & & S^{\gamma 1} \defn (S^a\backslash \{\gamma\}, \beta \seq S^s) \\ 
 & & S^\gamma \defn  \ben_{\Circle\alpha \in S^a}\mdl\Ep (S^a\backslash\{\mdl\alpha\},\alpha \seq \ ) \en 
   \ben_{\gamma=\mdl\alpha \imp \beta \in S^a}\Ep S^{\gamma 0} \en 
    (\mdl\Ap S^{\gamma 0}\imp \Ep S^{\gamma 1})  
\end{eqnarray*}
Then let 
\[\small
 \begin{array}{lll}
  \EpRSnot S & \defn & \left\{ 
                        \begin{array}{ll} 
                         S^\gamma & \text{if $S^s = \varnothing$} \\
                         \Ep (S^a \seq \, ) \en  S^\gamma & \text{if $S^s\neq                        \varnothing$} 
                        \end{array}
                      \right.                      
                         \\
  \ApRSnot S & \defn & \bof_{\gamma=\mdl\alpha \imp \beta \in S^a} 
   \mdl\Ap S^{\gamma 0}\en \Ap S^{\gamma 1}
 \end{array}
\]

\subsection{Soundness of the interpolant assignment}
 \label{secsoundness}

\begin{lemma}
For $\rsch=R\mdl$ all independent inductive properties and (DPP)$_\rsch$ hold. 
\end{lemma}
\begin{proof}
That (IPN)$_\rsch^\E$ and (IPN)$_\rsch^\A$ hold is clear. 
For (IPP)$_\rsch^\E$ and (IPP)$_\rsch^\A$, let $R$ be an instance 
\[
  \infer[=\ ]{S}{S_1}
  \infer{\Ga \seq \mdl\phi}{\Ga \seq \phi}
\]
of $\rsch=R\mdl$, $\EpR S$ and $\ApR S$ are defined as $\mdl\Ep S_1$ and $\mdl\Ap S_1$, respectively. Since $(\Ga \seq \Ep S_1)$ and 
$(\Ga, \Ap S_1 \seq \phi)$ belong to $\sqsRp$, an application of $\rsch$ shows that 
$\sqsRp \af \Ga \seq \EpR S$ holds, and an application of $\rsch$ followed by an application of $L\mdl$ shows that $\sqsRp \af \Ga,\ApR S \seq \mdl\phi$  holds. 

For (DPP)$_\rsch$, in a partition of $S$ in which $S^i$ is principal, $S^i$ has to be of the form $S^i=(\Ga^i \seq \mdl\phi)$, where $\Ga^i\subseteq \Ga$.  We have to show that $\sqsparRp \af S^r \cdot (\Ep S^i \seq \Ap S^i \mid \varnothing)$, which in this case means
\begin{equation} 
 \label{eqdpprmdl}
 \sqsparRp \af \Ga^r, \Ep S^i \seq \Ap S^i.
\end{equation}
Let $S_1^i=(\Ga^i \seq \phi)$. Thus $\Ga^r, \Ep S_1^i\seq \Ap S_1^i$ belongs to $\sqsparRp$. An application of $R\mdl$ followed by $L\mdl$ proves that 
$\sqsparRp \af \Ga^r, \mdl\Ep S_1^i \seq \mdl\Ap_1 S^i$. Since $S_1^i/S^i$ is an instance of $\rsch$, $\mdl\Ep S_1^i$ and $\mdl\Ap S_1^i$ are a conjunct and a disjunct of $\Ep S^i$ and $\Ap S^i$, respectively. This implies \eqref{eqdpprmdl}.
\end{proof}

\begin{lemma}
For $\rsch=L\mdl$ all independent inductive properties and (DPP)$_\rsch$ hold.
\end{lemma}
\begin{proof}
That (IPN)$_\rsch^\E$ and (IPN)$_\rsch^\A$ hold is clear. 
The proofs of (IPP)$^\E_\rsch$ and (IPP)$^\A_\rsch$ are similar to the one for $R\mdl$, but we include them nevertheless. 
For $R$ being an instance 
\[
  \infer[=\ ]{S}{S_1}
  \infer{\Ga,\mdl\psi \seq \mdl\phi}{\Ga,\psi \seq \mdl\phi}
\]
of $\rsch=L\mdl$, $\EpR S$ and $\ApR S$ are defined as $\mdl\Ep S_1$ and $\mdl\Ap S_1$, respectively. Thus $(\Ga,\psi \seq \Ep S_1)$ and 
$(\Ga, \psi, \Ap S_1 \seq \mdl\phi)$ belong to $\sqsRp$. An application of $R\mdl$ followed by an application of $\rsch$ to the former shows that $\sqsRp \af \Ga,\mdl\psi \seq \EpR S$ holds. Two applications of $\rsch$ to the latter show that $\sqsRp \af \Ga,\mdl\psi,\ApR S \seq \mdl\phi$ holds. 

For (DPP)$_\rsch$, we distinguish the case that $S^{is}$ is empty and that it is not. Note that $S^{ia}$ has to contain $\mdl\psi$ because $S^i$ is supposed to be principal for $L\mdl$. In the first case, $S^i=(\Ga^i,\mdl\psi \seq \ )$ for some $\Ga^i\subseteq \Ga$. We have to show that $\sqsparRp \af S^r \cdot (\Ep S^i \seq \Ap S^i \mid \varnothing)$, which in this case means
\begin{equation} 
 \label{eqdpplmdlone}
 \sqsparRp \af \Ga^r,\Ep S^i \seq \mdl\phi.
\end{equation}
Consider the partition of $S_1$ such that $S_1^i=(\Ga^i,\psi \seq \ )$. Thus 
$\Ga^r,\Ep S_1^i \seq \mdl\phi$ belongs to $\sqsparRp$ and an application of $L\mdl$ gives
$\sqsparRp \af \Ga^r,\mdl\Ep S_1^i \seq \mdl\phi$
Since $S_1^i/S^i$ is an instance of $\rsch$, $\mdl\Ep S_1^i$ is a conjunct of $\Ep S^i$. This implies \eqref{eqdpplmdlone}.

In the second case, $S^i=(\Ga^i,\mdl\psi \seq \mdl\phi )$ for some $\Ga^i\subseteq \Ga$. We have to show that $\sqsparRp \af S^r \cdot (\Ep S^i \seq \Ap S^i \mid \varnothing)$, which in this case means
\begin{equation} 
 \label{eqdpplmdltwo}
 \sqsparRp \af \Ga^r,\Ep S^i \seq \Ap S^i.
\end{equation}
Consider the partition of $S_1$ such that $S_1^i=(\Ga^i,\psi \seq \mdl\phi)$. Thus 
$\Ga^r,\Ep S_1^i \seq \Ap S_1^i$ belongs to $\sqsparRp$. An application of $R\mdl$ followed by $L\mdl$ proves $\sqsparRp \af \Ga^r,\mdl\Ep S_1^i \seq \mdl\Ap S_1^i$. 
Since $S_1^i/S^i$ is an instance of $\rsch$, $\mdl\Ep S_1^i$ and $\mdl\Ap S_1^i$ are a conjunct and a disjunct of $\Ep S^i$ and $\Ap S^i$, respectively. 
This implies \eqref{eqdpplmdltwo}. 
\end{proof}

\begin{lemma}
For $\rsch=R\mdl^\imp$ all independent inductive properties and (DPP)$_\rsch$ hold.
\end{lemma}
\begin{proof}
That (IPN)$_\rsch^\A$ holds is clear. 
For (IPP)$_\rsch^\E$ and (IPP)$_\rsch^\A$, let $R$ be an instance 
\[
  \infer[=\ ]{S}{S_1 & S_2}
  \infer{\Ga,\mdl\phi\imp \psi \seq \De}{\Ga \seq \phi & \Ga,\psi\seq \De}
\]
of $\rsch=R\mdl^\imp$. Since $\EpR S$ and $\ApR S$ are defined as $\Ep S_1 \en (\Ap S_1 \imp \Ep S_2)$ and $\Ap S_1 \en \Ap S_2$, respectively, we have to show that  
\[
 \sqsRp \af \Ga, \mdl\phi\imp \psi \seq \Ep S_1 \en (\Ap S_1 \imp \Ep S_2) \ \ \ \ 
 \sqsRp \af \Ga,\Ap S_1 \en \Ap S_2, \mdl\phi\imp \psi \seq \De.
\]
This follows from applications of $\rsch$ to sequents $(\Ga \seq \Ep S_1)$, $(\Ga,\psi\seq \Ep S_2)$, $(\Ga, \Ap S_1 \seq \phi)$ and $(\Ga, \psi, \Ap S_2 \seq \De)$, that all  belong to $\sqsRp$. 

For (IPN)$_\rsch^\E$ only the case that $\De$ is not empty is nontrivial. Thus suppose that $S$ is not principal for $R$ and that all sequents lower than $S$ satisfy the interpolant properties. Let $S_0=(S^a \seq \ )$. Hence $\EpRSnot S = \Ep S_0$. Because $S_0 \sml S$, $S_0$ satisfies the interpolant properties, and thus $S_0^a \seq \Ep S_0$ is derivable. Since $S^a=S_0^a$, it follows that $\af S^a \seq \EpRSnot S$. 
 
For (DPP)$_\rsch$, in a partition of $S$ in which $S^i$ is principal, $S^i$ has to be of the form $S^i=(\Ga^i,\mdl\phi\imp\psi \seq \De^i)$, where $\Ga^i\subseteq \Ga$ and $\De^i\subseteq \De$.  We treat the case that $\De^r$ is empty, the other case being analogous. We have to show that 
\begin{equation} 
 \label{eqdpprmdlimp}
\sqsparRp \af S^r \cdot (\Ep S^i \seq \Ap S^i).
\end{equation}
Let the partitions of the premisses be given by $S_1^i=(\Ga^i \seq \phi)$ and 
$S_2^i=(\Ga^i,\psi \seq \De^i)$. Thus $\Ga^r, \Ep S_1^i\seq \Ap S_1^i$ and 
$S_2^r \cdot (\Ep S_2^i\seq \Ap S_2^i)$ belong to $\sqsparRp$. 
As $(S_1^i \ S_2^i/S^i)$ is an instance of $\rsch$, $\Ep S^i$ contains a conjunct 
$\Ep S_1^i \en (\Ap S_1^i \imp \Ep S_2^i)$. Therefore $\sqsparRp$ derives $(\Ga^r,\Ep S^i \seq \Ap S_1^i)$ and $(\Ga^r,\Ep S^i \seq\Ap S_2^i)$. Since $\Ap S^i$ contains a disjunct $\Ap S_1^i \en \Ap S_2^i$, \eqref{eqdpprmdlimp} follows. 
\end{proof}

\begin{lemma}
For $\rsch=L\mdl^\imp$ all independent inductive properties and (DPP)$_\rsch$ hold.
\end{lemma}
\begin{proof}
For (IPP)$_\rsch^\E$ and (IPP)$_\rsch^\A$, let $R$ be an instance 
\[
  \infer[=\ ]{S}{S_1 & S_2}
  \infer{\Ga,\mdl\varchi,\mdl\phi\imp \psi \seq \De}{\Ga,\varchi \seq \mdl\phi & \Ga,\mdl\varchi,\psi\seq \De}
\]
of $\rsch=L\mdl^\imp$. Given the definition of $\EpR S$ and $\ApR S$ in this case, we have to show that  
\begin{eqnarray*}
 \sqsRp \af \Ga, \mdl\varchi,\mdl\phi\imp \psi \seq \mdl\Ep S_1 \en (\mdl\Ap S_1 \imp \Ep S_2) \\
 \sqsRp \af \Ga,\mdl\varchi,\mdl\Ap S_1, \Ap S_2, \mdl\phi\imp \psi \seq \De.
\end{eqnarray*}
This follows, via applications of $\rsch$, $R\mdl$ and $L\mdl$, from the fact that $(\Ga,\varchi \seq \Ep S_1)$ and $(\Ga,\mdl\varchi,\psi\seq \Ep S_2)$ as well as $(\Ga, \varchi, \Ap S_1 \seq \mdl\phi)$ and $(\Ga,\mdl\varchi, \psi, \Ap S_2 \seq \De)$ belong to $\sqsRp$. 

For (IPN)$_\rsch^\E$ we have to show three things, under the assumption that all sequents lower than $S$ satisfy the interpolant properties:
\begin{enumerate}
\item[(1)]
$\af S^a \seq \Ep (S^a \seq\ )$ in case $S^s \neq \varnothing$;
\item[(2)] 
$\af S^a \seq \mdl \Ep (S^a\backslash \{\mdl\alpha\},\alpha\seq \ )$ for all $\mdl\alpha \in S^a$;
\item[(3)] 
$\af S^a \seq \Ep S^{\gamma 0}\en (\mdl \Ap S^{\gamma 0}\imp \Ep S^{\gamma 1})$ for all $\gamma=\mdl\alpha\imp\beta$ in $S^a$. 
\end{enumerate}
(1) follows immediately from the assumption that all sequents lower than $S$ satisfy the interpolant properties and the fact that $(S^a \seq\ )\sml S$, which holds if $S^s\neq\varnothing$. 

(2) Note that $\big(S^a\backslash \{\mdl\alpha\},\alpha \seq \Ep (S^a\backslash \{\mdl\alpha\},\alpha\seq\ )\big)$ is derivable because all sequents lower than $S$ are assumed to satisfy the interpolant properties. An application of $R\mdl$ followed by an application of $L\mdl$ gives the desired result. 

(3) Assume that $S=(\Pi,\gamma \seq \Sig)$, where $\gamma=\mdl\alpha\imp\beta$. We use the following abbreviations:
\begin{eqnarray*}
 S_\alpha & \defn & 
  S^{\gamma 0} =(S^a\backslash \{\gamma\} \seq \mdl\alpha) = (\Pi \seq \mdl\alpha) \\
 S_\beta & \defn & 
  S^{\gamma 1}=(S^a\backslash \{\gamma\},\beta \seq S^s)=(\Pi,\beta \seq \Sig).
\end{eqnarray*}
Thus we have to show that $(S^a \seq \Ep S_{\alpha})$ and 
$(S^a,\mdl\Ap S_{\alpha} \seq \Ep S_{\beta})$ are derivable. The former follows immediately from the fact that $S_\alpha$ and $S_\beta$ satisfy the interpolant properties, since they are below $S$. By the same fact, $(S_{\beta}^a \seq \Ep S_{\beta})$ is derivable. As $(S_{\beta}^a \seq \Ep S_{\beta}\big)$ is equal to $(\Pi,\beta \seq \Ep S_{\beta})$, sequent $(\Pi,\beta, \mdl \Ap S_{\alpha} \seq \Ep S_{\beta})$ is derivable too. An application of $L\mdl^\imp$ to this sequent and 
$(S_{\alpha}^a,\Ap S_{\alpha} \seq S_{\alpha}^s) =(\Pi, \Ap S_{\alpha} \seq \mdl\alpha)$, which by the same fact is derivable, shows that $\af \Pi,\gamma, \mdl\Ap S_{\alpha} \seq \Ep S_{\beta}$, which is what we had to show. 

For property (IPN)$_\rsch^\A$ we use the same notation as in (3) of the previous case and have to show that 
\[
 \af S^a, \mdl\Ap S_\alpha, \Ap S_\beta \seq S^s.  
\]
By the assumption that all sequents lower than $S$ satisfy the interpolant properties, both $(S_\alpha^a,\Ap S_\alpha \seq S_\alpha^s) = (\Pi,\Ap S_\alpha \seq \mdl\alpha)$ and $(S_\beta^a,\Ap S_\beta \seq S_\beta^s)=(\Pi,\beta,\Ap S_\beta \seq \Sig)$ are derivable. Clearly, applications of weakening, $L\mdl$ and $L\mdl^\imp$ show that 
$\af \Pi,\mdl\alpha\imp\beta, \mdl\Ap S_\alpha, \Ap S_\beta \seq \Sig$, which is what we had to show. 

For (DPP)$_\rsch$, in a partition of $S$ in which $S^i$ is principal, $S^i$ has to be of the form $S^i=(\Ga^i,\mdl\varchi, \mdl\phi\imp\psi \seq \De^i)$, where $\Ga^i\subseteq \Ga$ and $\De^i\subseteq \De$.  We have to show that 
\begin{equation} 
 \label{eqdpplmdlimp}
\sqsparRp \af S^r \cdot (\Ep S^i \seq \Ap S^i \mid \varnothing)
\end{equation}
We treat the case that $\De^r=\varnothing$, the other case is analogous. 
Let the partitions of the premisses be given by $S_1^i=(\Ga^i,\varchi \seq \phi)$ and 
$S_2^i=(\Ga^i,\mdl\varchi,\psi \seq \De^i)$. Thus $(\Ga^r, \Ep S_1^i\seq \Ap S_1^i)$ and 
$S_2^r \cdot (\Ep S_2^i\seq \Ap S_2^i)$ belong to $\sqsparRp$. Hence $\sqsparRp$ derives $(\Ga^r, \mdl\Ep S_1^i\seq \mdl\Ap S_1^i)$. 
As $(S_1^i\ S_2^i/S^i)$ is an instance of $\rsch$, $\Ep S^i$ contains a conjunct 
$\mdl\Ep S_1^i \en (\mdl\Ap S_1^i \imp \Ep S_2^i)$. Therefore $\sqsparRp$ derives $(\Ga^r,\Ep S^i\seq \mdl\Ap S_1^i)$ and $(\Ga^r,\Ep S^i\seq \Ap S_2^i)$. Since $\Ap S^i$ contains a disjunct $\mdl\Ap S_1^i \en \Ap S_2^i$, \eqref{eqdpplmdlimp} follows. 
\end{proof}

\begin{lemma}
For all modal and implication-modal rules $\rsch$, the interpolant assignment satisfies (DPN)$_\rsch$. 
\end{lemma}
\begin{proof}
(DPN)$_\rsch$ is proved with induction to the depth of the derivation $\cald$ of $S$ for all rules at the same time. Thus we have to show for any partition $(S^r,S^i)$ of $S$ such   that $S^i$ is nonprincipal for the last inference $R$ of $\cald$: if all sequents below $S$ satisfy the interpolant properties, then 
\begin{equation}
 \label{eqdpn}
 \af S^r \cdot (\Ep S^i \seq \Ap S^i \mid \varnothing).
\end{equation} 

Note that for the rules in \DY\, (DPN)$_\rsch$ has already been established in \citep{iemhoff17}. If $\cald$ has depth one, $S$ is an instance of an axiom. If the axiom is $L\bot$, then because $S^i$ is not an instance of it, $\bot \in S^{ra}$. Hence \eqref{eqdpn}. If the axiom is {\it Ax}, say $S = (\Ga,q\seq q,\De)$, then if $q \in S^{ra} \cap S^{rs}$, we have \eqref{eqdpn}. Therefore assume otherwise: $q \not\in S^{ra} \cap S^{rs}$.
Because $S^i$ is not an instance of {\it Ax}, $q\not\in S^{ia} \cap S^{is}$. Thus $q\in S^{ia} \cap S^{rs}$ or $q \in S^{ra} \cap S^{is}$. Since $q$ occurs in $S^r$, $q \neq p$. Therefore $\af q \seq S^{rs}$ and $q \in S^{ia}$, and thus $\af \Epat S^i \seq q$, or $\af S^{ra} \seq q$ and 
$q \in S^{is}$, and thus $\af q \seq \Apat S^i$. In both cases \eqref{eqdpn} follows, which completes the basis case of the induction. 

Suppose the derivation of $S$ has depth bigger than one and ends with an application of $\rsch$. We show that (DPN)$_\rsch$ holds for $S$. If $\rsch$ is not a modal or implication-modal rule, then we already know that (DPN)$_\rsch$ holds for all $S$. Therefore suppose $\rsch$ is such a rule. We consider the four rules separately. 
 
If $\rsch$ is $R\mdl$, $S = (\Ga \seq \mdl\phi)$ and the premise of the last inference is $S_1 =(\Ga\seq\phi)$. Since we consider (DPN)$_\rsch$, $S^i=(\Ga^i\seq \ )$. We have to show that $\sqsparRp \af \Ga^r, \Ep S^i \seq \mdl\phi$. Consider the partition of $S_1$ with $S_1^i = S^i$. We have $\sqsparRp \af \Ga^r, \Ep S_1^i \seq \phi$ by the induction hypothesis, which implies what we had to show via the consecutive application of $R\mdl$ and $L\mdl$.

If $\rsch$ is $L\mdl$, $S = (\Ga,\mdl\psi \seq \mdl\phi)$ and the premise of the last inference is $S_1 =(\Ga,\psi \seq\mdl\phi)$. Consider a partition of $S$ such that $S^i$ is not principal for $L\mdl$. Thus $S^i=(\Ga^i\seq \mdl\phi)$ or $S^i=(\Ga^i\seq \ )$ for some $\Ga^i\subseteq \Ga$. We have to show that $\sqsparRp \af \Ga^r, \mdl\psi, \Ep S^i \seq \Ap S^i$ in the first case, and $\sqsparRp \af \Ga^r, \mdl\psi,\Ep S^i \seq \mdl\phi$ in the second case. In both cases we can partition $S_1^i$ such that $S_1^i=S^i$. If $S^i=(\Ga^i\seq \mdl\phi )$, we have $\sqsparRp \af \Ga^r,\psi,\Ep S_1^i \seq \mdl\phi$ by the induction hypothesis. An application of $L\mdl$ gives 
$\sqsparRp \af \Ga^r, \mdl\psi,\Ep S_1^i \seq \mdl\phi$, and since $S_1^i=S^i$, this is what we had to show. The case that $S^i=(\Ga^i\seq \ )$ is analogous. 

If $\rsch$ is $R\mdl^\imp$, the last inference is of the following form:
\[
  \infer[=\ ]{S}{S_1 & S_2}
  \infer{\Ga,\mdl\phi \imp \psi \seq \De}{\Ga\seq \phi & \Ga,\psi \seq \De}
\]
As $\rsch$ is not applicable to $S^i$, $S^i = (\Ga^i \seq \De^i)$ and 
$S^r =(\Ga^r,\mdl\phi \imp \psi \seq \De^r)$. We treat the case that $\Delta^r=\varnothing$, the other case is analogous. Thus to show \eqref{eqdpn}, we have to show that 
\begin{equation}
 \label{eqdpnrim}
 \af \Ga^r,\mdl\phi\imp\psi,\Ep S^i \seq \Ap S^i.
\end{equation}
Let $S_1^i=(\Ga^i \seq \ )$ and $S_2^i=S^i$. By the induction hypothesis, 
\[
 \af \Ga^r,\Ep S_1^i \seq \phi \ \ \ \ 
 \af \Ga^r,\psi,\Ep S_2^i \seq \Ap S_2^i.
\]
Therefore $\af \Ga^r,\mdl\phi\imp\psi,\Ep S_1^i,\Ep S_2^i \seq \Ap S_2^i$. Since $S_2^i=S^i$, for \eqref{eqdpnrim} it suffices to show that $\Ep S^i$ derives $\Ep S_1^i$. 
If $S^{is}=\varnothing$, then $S^i=S_1^i$, and if $S^{is}\neq \varnothing$, then $\Ep S_1^i$ is a conjunct of $\EpRSnot S^i$, which is a conjunct of $\Ep S^i$. 

If $\rsch$ is $L\mdl^\imp$, the last inference is of the following form:
\[
  \infer[=\ ]{S}{S_1 & S_2}
  \infer{\Ga,\mdl\varchi,\mdl\phi \imp \psi \seq \De}{\Ga, \varchi \seq \mdl\phi & \Ga,\mdl\varchi,\psi \seq \De}
\]
We distinguish three cases. 

(1) If $S^i = (\Ga^i \seq \De^i)$, then let $S_1^i=(\Ga^i \seq \ )$ and $S_2^i=S^i$. We treat the case that $\Delta^r\neq\varnothing$, the other case is analogous. Thus to show \eqref{eqdpn}, we have to show that 
\begin{equation}
 \label{eqdpnlimone}
 \af \Ga^r,\mdl\varchi, \mdl\phi\imp\psi,\Ep S^i \seq \De^r.
\end{equation}
By the induction hypothesis, 
\[
 \af \Ga^r,\varchi, \Ep S_1^i \seq \mdl\phi \ \ \ \ 
 \af \Ga^r,\mdl\varchi, \psi,\Ep S_2^i \seq \De^r.
\]
Therefore $\af \Ga^r,\mdl\varchi, \mdl\phi\imp\psi,\Ep S_1^i,\Ep S_2^i \seq \De^r$. If $\De^i=\varnothing$, then $S^i=S_1^i=S_2^i$, which implies  \eqref{eqdpn}. Therefore assume that $\De^i = S^{is}$ is not empty. Thus $\Ep S_1^i$ is a conjunct of 
$\EpRSnot S^i$, which is a conjunct of $\Ep S^i$. Together with $S_2^i=S^i$ this implies \eqref{eqdpnlimone}.

(2) If $S^i = (\Ga^i,\mdl\varchi \seq \De^i)$, then let $S_1^i=(\Ga^i,\varchi \seq \ )$ and $S_2^i=S^i$. We treat the case that $\Delta^r=\varnothing$, the other case is analogous. Thus to show \eqref{eqdpn}, we have to show that 
\begin{equation}
 \label{eqdpnlimtwo}
 \af \Ga^r,\mdl\phi\imp\psi,\Ep S^i \seq \Ap S^i.
\end{equation}
By the induction hypothesis, $\af \Ga^r, \Ep S_1^i \seq \mdl\phi$ and $\af \Ga^r, \psi,\Ep S_2^i \seq \Ap S_2^i$ are derivable. An application of $L\mdl^\imp$ gives 
\[
 \af \Ga^r,\mdl\phi\imp \psi,\mdl\Ep S_1^i,\Ep S_2^i \seq \Ap S_2^i.
\] 
Note that because $S_1^i=(S^i\backslash \{\mdl\varchi\},\varchi \seq \ )$, $\mdl\Ep S_1^i$ is a conjunct of $\EpRSnot S^i$, which is a conjunct of $\Ep S^i$. 
Since also $\Ep S_1^i\seq\mdl\Ep S_1^i$ is derivable and $S_2^i=S^i$, we obtain  \eqref{eqdpnlimtwo}. 

(3) If $S^i = (\Ga^i,\mdl\phi \imp \psi \seq \De^i)$, then let $S_1^i=(\Ga^i \seq \mdl\phi)$ and $S_2^i=(\Ga^i,\psi \seq \De^i)$. We treat the case that $\Delta^r=\varnothing$, the other case is analogous. Thus to show \eqref{eqdpn}, we have to show that 
\begin{equation}
 \label{eqdpnlimthree}
 \af \Ga^r,\mdl\varchi,\Ep S^i \seq \Ap S^i.
\end{equation}
By the induction hypothesis, both $(\Ga^r,\varchi, \Ep S_1^i \seq \Ap S_1^i)$ and $S_3=(\Ga^r,\mdl\varchi, \Ep S_2^i \seq \Ap S_2^i)$ are derivable. 
Thus so are $S_4=(\Ga^r,\varchi, \Ep S_1^i \seq \mdl\Ap S_1^i)$ and $(\Ga^r,\mdl\varchi, \Ep S_1^i \seq \mdl\Ap S_1^i)$ 
by applications of $R\mdl$ and $L\mdl$. 
An application of $L\mdl^\imp$ to $S_4$ and $S_3$ shows the derivability of $(\Ga^r,\mdl\varchi,\Ep S_1^i, \mdl\Ap S_1^i\imp\Ep S_2^i \seq \Ap S_2^i)$. Therefore 
\[
 \af \Ga^r,\mdl\varchi, \Ep S_1^i, \mdl\Ap S_1^i\imp\Ep S_2^i \seq 
 \mdl\Ap S_1^i \en \Ap S_2^i.
\] 
Let $\gamma=\mdl\phi\imp\psi$. Since $S_1^i=(S^i)^{\gamma 0}$ and $S_2^i=(S^i)^{\gamma 1}$, $\Ep S_1^i$ and $\mdl\Ap S_1^i\imp\Ep S_2^i$ are conjuncts of $\EpRSnot S^i$, which is a conjunct of $\Ep S^i$.   And since $\mdl\Ap S_1^i \en \Ap S_2^i$ is a disjunct of $\ApRSnot S^i$ and thus of $\Ap S^i$, we can conclude \eqref{eqdpnlimthree}.  
\end{proof}

\subsection{Uniform interpolation for Lax Logic}

\textcolor{cyan}{The proof of the following theorem is not correct, as it is based on Corollary 1, which proof is not correct, see the discussion below Corollary 1.} 

\begin{theorem}
\PLL\ has uniform interpolation. 
\end{theorem}
\begin{proof}
As the results in Section~\ref{secsoundness} show, the interpolant assignment for $\GLL$ given in Section~\ref{secassignment} is sound. Theorem~\ref{thmsufficientmod} then implies the result. 
\end{proof}

\section{Conclusion}
In this paper we have developed proof theory for propositional Lax Logic \PLL. In particular, we have provided a cut free terminating sequent calculus \DYLL\ for the logic, which is used to prove that Lax Logic has uniform interpolation. The latter was established by a proof-theoretic method that defines uniform interpolants in terms of the rules of the calculus. 

There are many open questions related to uniform interpolation in intuitionistic modal logics. We hope that the methods in this paper and the closely related \citep{iemhoff17,jalali&tabatabai2018b} will be useful in the study of other logics in this class, and perhaps also in the investigation of some of the more general questions in Universal Proof Theory.

\paragraph{Acknowledgements}
We thank an anonymous referee for the constructive comments on an earlier draft of this paper and Iris van der Giessen for several useful conversations on Lax Logic, as well as Raheleh Jalali and Amir Tabatabai for numerous discussions about proof theory and uniform interpolation.

\end{document}